\documentclass[12pt]{amsart}

\makeatletter

\renewcommand{\subsection}{\@startsection{subsection}{2}%
  \z@{.5\linespacing\@plus.7\linespacing}{-.5em}%
  {\normalfont\bfseries\S\,}}

\makeatother

\usepackage{amsmath}
\usepackage{amssymb}
\usepackage{enumerate}
\usepackage{euscript}
\usepackage{amscd}
\usepackage[all]{xy}

\newtheorem{thm}{Theorem}[section]
\newtheorem{lem}[thm]{Lemma}
\newtheorem{prop}[thm]{Proposition}
\newtheorem{cor}[thm]{Corollary}
\theoremstyle{definition}
\newtheorem{defn}[thm]{Definition}
\newtheorem{rem}[thm]{Remark}

\numberwithin{equation}{section}
\theoremstyle{remark}

\newcommand{\bbc}{{\mathbb C}}

\newcommand{\bbq}{{\mathbb Q}}
\newcommand{\bbr}{{\mathbb R}}

\newcommand{\bbz}{{\mathbb Z}}

\newcommand{\del}{{\delta}}

\newcommand{\gC}{{\mathfrak C}}

\newcommand{\gii}{{\mathfrak e}}

\newcommand{\gM}{{\mathfrak M}}

\newcommand{\gp}{{\mathfrak p}}

\newcommand{\gr}{{\mathfrak r}}

\newcommand{\cK}{{\mathcal K}}

\newcommand{\cN}{{\mathcal N}}

\newcommand{\co}{{\mathcal O}}
\newcommand{\cD}{{\mathcal D}}

\newcommand{\aff}{{\operatorname {Aff}}}

\newcommand{\n}{{\operatorname {N}}}

\newcommand{\gl}{{\operatorname{GL}}}

\newcommand{\spl}{{\operatorname{SL}}}

\newcommand{\sst}{{\operatorname{ss}}}

\newcommand{\vol}{{\operatorname{vol}}}
\newcommand{\ord}{{\operatorname{ord}}}
\newcommand{\sym}{{\operatorname{Sym}}}

\newcommand{\Z}{\bbz}
\newcommand{\Q}{\bbq}
\newcommand{\R}{\bbr}
\newcommand{\C}{\bbc}

\newcommand{\gMf}{\gM_{\rm f}}
\newcommand{\gMi}{\gM_{\infty}}
\newcommand{\gMc}{\gM_\C}
\newcommand{\gMr}{\gM_\R}
\newcommand{\gMdy}{\gM_{\text {dy}}}

\newcommand{\nr}{\EuScript N}

\newcommand{\twtw}[4]{\left(\begin{array}{cc}{#1}&{#2}\\{#3}&{#4}\\\end{array}\right)}

\begin{document}

\title[proportional constants]
{On proportional constants of the mean value
of class numbers of quadratic extensions}
\author[Takashi Taniguchi]{Takashi Taniguchi}
\address{Graduate School of Mathematical Sciences\\ University of Tokyo\\
3--8--1 Komaba Megoro-ku\\ Tokyo 153-0041\\ JAPAN}
\email{tani@ms.u-tokyo.ac.jp}
\date{\today}

\begin{abstract}
In this article, we slightly refine the mean value theorem
for the class number of quadratic extensions obtained by
Goldfeld-Hoffstein and Datskovsky.
We determine all the proportional constants of the mean value
with respect to the local conditions including dyadic places.
\end{abstract}

\maketitle

\section{Introduction}\label{sec:int}

For a number field $k$ let $h_k$ 
and $R_k$ be the class number and the regulator,
respectively.
For a finite extension $F/k$ of number fields,
let $\cN(\Delta_{F/k})$ be the absolute norm of the relative
discriminant.
In this article, we will give a refinement of
Datskovsky's \cite{dats} mean value theorem
of $h_FR_F$ with respect to $\cN(\Delta_{F/k})$
for certain families of quadratic extensions $F$
of a fixed number field $k$.
Our main theorem is Theorem \ref{maintheorem}.

The roots of this topic traces to Gauss.
Let $h_d$ be the number of $\spl(2,\Z)$-equivalence classes
of primitive integral binary quadratic forms
which are either positive definite or indefinite.
In \cite[Sections 302, 304]{gauss}
Gauss gives conjectures for the asymptotic property of
the average number of $h_d$.
This conjecture was first proved
by Lipschitz for the imaginary case,
and by Siegel for the real case.
Siegel \cite{siegele} also proved
a density theorem for integral equivalence classes
of quadratic forms in general.

M. Sato and T. Shintani formulated this kind of density problems
using the notion of prehomogeneous vector spaces.
In \cite{Shintanib} Shintani considered
the zeta functions associated with
the space of quadratic forms.
There he reproved the density theorem
for the case of binary quadratic forms
and improved the error estimate.
However, in Gauss' conjecture
all integers $d$ are allowed,
and if $d=m^2d'$ and $d'$ is square free integer,
$h_d$ and $h_{d'}$ are related by a simple relation.
Therefore we are counting essentially the same object infinitely many times.
This ambiguity was first removed by Goldfeld-Hoffstein.
In \cite{goho}, they give the mean value of $h_kR_k$ of quadratic fields $k$
by using Eisenstein series of half-integral weight.

Later Datskovsky \cite{dats} investigated this subject
from the view of prehomogeneous vector spaces.
He constructed the necessary local theory,
and combined it with Shintani's global theory \cite{Shintanib}.
He then generalized the Goldfeld-Hoffstein mean value theorem
to quadratic extensions of an arbitrary fixed number field.
In the proof \cite{dats}, he also showed the existence of a
the mean value of the family of quadratic extensions
which have a given local behavior
at a finite number of places of $k$,
and computed the proportional constants to the whole mean value
explicitly for many type of local conditions.

However, the computation of the proportional constants with respect to
quadratic ramified extensions have some difficulty
especially in the case of wild ramification at dyadic places
(those dividing the place of $\Q$ at 2),
and for certain types he gave the sum of proportional constants
over certain arithmetically similar conditions instead.
The purpose of this article is to determine
each unevaluated proportional constant
to complete Datskovsky's work.

We follow Datskovsky's approach.
Our method is a natural modification of Kable-Yukie \cite{kayuII},
in which another prehomogeneous vector space is handled
and a new mean value theorem is obtained.
In \cite{tanic} we find a mean value of $h_F^2R_F^2$
with respect to $\cN(\Delta_{F/k})$
for certain families of quadratic extensions $F$
of a fixed number field $k$.
We can also apply the computation
in this article to that of \cite{tanic}.

The contents of this article is as follows.
In Section \ref{sec:rev}
we recall the definition of the space of binary quadratic forms
and review its fundamental properties.
For later purposes, we define the representation over any ring.
In Section \ref{sec:mea}
we first state the main theorem of this article,
which is a refinement of Datskovsky's mean value theorem \cite{dats}.
After that we discuss the relation
between the space of binary quadratic forms
and the main theorem,
and point out what is the remaining task.
This reduces to determine certain orbital volumes
in a vector space over a non-archimedean local field.
The explicit computation of the volume
is carried out in the final section.

\bigskip

\noindent
{\bf Notations.}
For a finite set $X$ we denote by $\# X$ its cardinality.
The standard symbols $\Q$, $\R$, $\C$ and $\Z$ will denote respectively
the rational, real and complex numbers and the rational integers.
If $V$ is a scheme defined over a ring $R$ and $S$ is an $R$-algebra
then $V_S$ denotes its $S$-rational points.
For a finite field extension $L/F$,
let $\nr_{L/F}$ denote the norm map.

\section{Review of the space of binary quadratic forms}\label{sec:rev}
In this section, we briefly review
the space of binary quadratic forms.
Let $G=\gl(1)\times\gl(2)$ and $V=\sym^2\aff^2$.
We consider the natural action of $\gl(2)$ on $V$.
We define the action of $\gl(1)$ on $V$ by the usual
scalar multiplication.
This defines a representation of $G$ on $V$.
We regard $V$ as the space of
binary quadratic forms of variables $v=(v_1,v_2)$.
Elements of $V$ are expressed in the form
\begin{equation*}
x=x(v)=x(v_1,v_2)=x_0v_1^2+x_1v_1v_2+x_2v_2^2.
\end{equation*}
For the rest of this paper, we express elements of $G$
as
\begin{equation*}
g=(t,g_2),\qquad t\in\gl(1),
\quad g_2=\twtw {g_{211}}{g_{212}}{g_{221}}{g_{222}}\in\gl(2).
\end{equation*}
Then the action of $g$ is given by $x(v)\mapsto tx(vg_2)$,
regarding $v$ as a row vector.

Let $P(x)=x_1^2-4x_0x_2$, the discriminant of $x(v)$.
Let $\chi_i\ (i=1,2)$ be characters on $G$ defined by
$\chi_1(g)=t, \chi_2(g)=\det(g_2)$, respectively.
We define the character $\chi$ on $G$ by $\chi=\chi_1\chi_2$.
Then $P$ becomes a relative invariant polynomial
with respect to the character $\chi^2$
i.e. $P(gx)=\chi(g)^2P(x)$ for all $g\in G, x\in V$.
If $x$ is a rational point of $V$,
we can consider the stabilizer of $x$
as a group scheme,
in the sense of \cite[p.3]{mufo}.
We denote this group scheme by $G_x$.

Let $x=v_1^2+a_1v_1v_2+a_2v_2^2$ be a rational point of $V$.
For these forms of $x$, we put
\begin{equation*}
A_x(c,d)=\twtw cd{-a_2d}{c+a_1d}
\end{equation*}
and then define the subgroup $N_x$ of $G$ by
\begin{equation*}
N_x=\left\{n_x(c,d)=(\det(A_x(c,d))^{-1},A_x(c,d))
\mid A_x(c,d)\in\gl(2)\right\}.
\end{equation*}
By computation, one could see that
$N_x$ is a subgroup of $G_x$.

We now consider the representation $(G,V)$
over a fixed field $K$.
Let $V^\sst=\{x\in V\mid P(x)\not=0\}$
and call the set of semi-stable points.
By definition,
$x\in V_K$ is an element of $V_K^\sst$
if and only if $x(v)$ has distinct roots in $\mathbb P^1_{\overline K}$.

\begin{defn}
For $x\in V_K^\sst$, we define
\begin{align*}
Z_x		&=\mathrm{Proj}\, K[v_1,v_2]/(x(v)),\\
K(x)	&=\Gamma(Z_x,\co_{Z_x}).
\end{align*}
\end{defn}
Note that $K(x)$ may not be a field.
Since $V_K^\sst$ is the set of $x$ such that
$F_x$ does not have a multiple root,
$Z_x$ is a reduced scheme over $k$ and
$K(x)$ is a separable quadratic algebra over $K$.
Also note that if $x(v)$ is irreducible over $K$,
then $K(x)$ is the splitting field of $x(v)$.
The following proposition is easy to see and we omit the proof.
\begin{prop}\label{rod}
The map $x\mapsto K(x)$ gives a bijection
between the set of rational orbits $G_k\backslash V_K^\sst$
and the set of equivalence class of separable quadratic algebra over $K$.
\end{prop}

\begin{rem}
For $x\in V_K^\sst$, the identity component of the stabilizer of $x$,
denoted by $G_x^\circ$,
is isomorphic to $R_{K(x)/K}(\mathbb G_m)$,
the scalar restriction of $\mathbb G_m$
with respect to $K(x)/K$.
Hence if $K$ is a number field and $K(x)$ is a quadratic extension of $K$,
then the unnormalized Tamagawa number of $G_x^\circ$ is
the constant multiple of $h_{K(x)}R_{K(x)}$,
and this is the reason why the study of representation $(G,V)$
leads the mean value theorem of class number times regulator
of quadratic extensions.
\end{rem}

The proof of the following proposition is
straightforward and we omit the detail.

\begin{prop}\label{Nstr}
Let $x=v_1^2+x_1v_1v_2+x_2v_2^2\in V_K^\sst$
be an irreducible form over $k$
and $x=(v_1+\alpha_1 v_2)(v_1+\alpha_2v_2)$
for $\alpha_1,\alpha_2\in K(x)$.
Then the map
\begin{equation*}
\psi_{x,K}\colon
N_{x\, K}\longrightarrow K(x)^\times\quad\text{as}\quad
n_x(c,d)\longmapsto c+\alpha_1d
\end{equation*}
gives an isomorphism of these groups.
Moreover the following diagram is commutative:
\begin{equation*}
\begin{CD}
N_{x\, K}			@>{\chi_2}>>	K^\times	\\
@V{\psi_{x,K}}VV		@|				\\
K(x)^\times			@>{\nr_{K(x)/K}}>>	K^\times	\\
\end{CD}
\end{equation*}
\end{prop}
Note that in the case of the proposition above,
$N_x$ coincides with the identity component of the stabilizer $G_x$
as an algebraic group over $K$.

\section{The mean value theorem}\label{sec:mea}
In this section, we state the main theorem of this article,
which is a slight refinement of Datskovsky's \cite{dats}.
After that, we point out what is the remaining task.
The explicit computations are done in the next section.

To state the main theorem, we prepare some notations.
For the rest of this article we fix a number field $k$.
Let $\gM$,
$\gM_{\infty}$, $\gM_{\text{f}}$,
$\gM_{\R}$ and $\gM_{\C}$
denote respectively the set of all places of $k$, all infinite
places, all finite places,
all real places and all complex places.
We let $r_1$, $r_2$, $h_k$, $R_k$ and $e_k$ be
respectively the number of real places, the number of complex
places, the class number, the regulator and the number of roots of
unity contained in $k$. We set
$\gC_k=2^{r_1}(2\pi)^{r_2}h_kR_ke_k^{-1}$.
We write by $\zeta_k(s)$ the Dedekind zeta function of $k$.
For $v\in\gM$ let $k_v$ denote the completion of $k$ at $v$ and
$|\ |_v$ the normalized absolute value of $k_v$.
For $v\in\gMf$, we denote by $q_v$
the cardinality of the residue field at $v$.
For a finite extension $L_v/k_v$ at $v\in\gMf$,
let $\Delta_{L_v/k_v}$ denote the relative discriminant.

Let $S\supset\gMi$ be a finite set of places.
We consider $S$-tuples $L_S=(L_v)_{v\in S}$
where each $L_v$ is a separable quadratic algebra of $k_v$.
Let $F$ be a quadratic extension of $k$.
We write $F\approx L_v$ to mean that
$F\otimes_k k_v$ is isomorphic to $L_v$ as a $k_v$-algebra.
For  $L_S=(L_v)_{v\in S}$ 
we shall write $F\approx L_S$ if and only if
$F\approx L_v$ for all $v\in S$.
For $v\in\gM$ and $L_v$ a separable quadratic algebra over $k_v$,
we define the constant $\gii_v(L_v)$ as follows.
\begin{defn}\label{gii}
\begin{enumerate}[{\rm (1)}]
\item
If $v\in\gMf$, then we  put
\begin{equation*}
\gii_v(L_v)
=
\begin{cases}
	2^{-1}(1-q_v^{-2})
		&	L_v\cong k_v\times k_v,\\
	2^{-1}(1-q_v^{-1})^2
		&	L_v\text{ is quadratic unramified},\\
	2^{-1}|\Delta_{L_v/k_v}|_v^{-1}(1-q_v^{-1})(1-q_v^{-2})
		&	L_v\text{ is quadratic ramified}.\\
\end{cases}
\end{equation*}
\item
For $v\in\gMr$, we define $\gii_v(L_v)=1/4$ if $L_v\cong\R\times\R$,
and $\gii_v(L_v)=1/(2\pi)$ if $L_v\cong\C$.
\item
For $v\in\gMc$, we define $\gii_v(L_v)=1/(4\pi^2)$.
\end{enumerate}
\end{defn}
Let $\cN(\Delta_{F/k})$ denote the absolute norm
of the relative discriminant $\Delta_{F/k}$ of $F/k$.
Then the following, a refinement of \cite{dats},
is a main result of this paper.
\begin{thm}\label{maintheorem}
Let $S\supset\gMi$ and $L_S=(L_v)_{v\in S}$ is a $S$-tuple.
Then we have
\begin{equation*}
\lim_{X\to\infty}\frac{1}{X^{3/2}}
\sum_{\underset{\cN(\Delta_{F/k})\leq M}{[F:k]=2,\ F\approx L_S}}
h_FR_F
=	\frac{e_k\gC_k^2\zeta_k(2)}{3\cdot2^{r_1+r_2-1}}
		\prod_{v\in S}\gii_v(L_v)
		\prod_{v\in\gM\setminus S}(1-q_v^{-2}-q_v^{-3}+q_v^{-4}).
\end{equation*}
\end{thm}
The contribution of this article is that
we determine the constant $\gii_v(L_v)$
for $v\in\gMf$ and $L_v$ quadratic ramified extension solely,
whereas in \cite{dats}, the sum of $\gii_v(L_v)$
for $L_v$'s with the same relative discriminants were calculated.

We will discuss the relation between
the prehomogeneous vector space $(G,V)$
and the constant $\gii_v(L_v)$.
For the rest of this paper we assume $v\in\gMf$.
We denote by $\co_v$ the integer ring of $k_v$.
We consider the prehomogeneous vector space $(G,V)$ over $\co_v$
and over $k_v$.
Let $\cK_v=G_{\co_v}$,
which is the standard maximal compact subgroup of $G_{k_v}$.
On $V_{k_v}$ we consider the additive Haar measure under which
$\vol(V_{\co_v})=1$.

As in Proposition \ref{rod}, the set of rational orbits
$G_{k_v}\backslash V_{k_v}^\sst$ corresponds bijectively
to the set of separable quadratic algebra of $k_v$.
Following \cite{dats} we select and fix
a representative element from each orbit.
\begin{defn}
For each of $G_{k_v}$-orbits in $V_{k_v}^\sst$,
we choose and fix an element $x$ as follows.
\begin{enumerate}[(1)]
\item
For the orbit corresponding to the algebra $k_v\times k_v$,
we set $x=v_1v_2$.
\item
For any orbit corresponding to a quadratic extension $L$, 
set $x=\nr_{L/k_v}(v_1+\theta v_2)$ where
$\theta$ is a generator of the integer ring $\co_L$ of $L$ over $\co_v$.
\end{enumerate}
We call such fixed orbital representatives as
the {\em standard orbital representatives}.
\end{defn}
Note that if $x\in V_{k_v}^\sst$ is a standard representative,
then $P(x)$ generates the relative discriminant $\Delta_{k_v(x)/k_v}$.

Let $x\in V_{k_v}^\sst$ be any standard representative.
Then it is shown in \cite{dats} that
$\gii_v(k(x))=\vol(\cK_vx)$.
He then computed $\vol(\cK_vx)$
if $x$ corresponds to $k_v\times k_v$ or
to the quadratic unramified extension,
and $\sum_{x}\vol(\cK_vx)$ where $x$ runs through 
all standard representatives
corresponding to quadratic ramified extensions
with 
the given relative discriminant
$\Delta_{k_v(x)/k_v}$.
Hence, to obtain Theorem \ref{maintheorem},
we will determine the value $\vol(\cK_vx)$ for
each standard representative $x$
corresponding to a quadratic ramified extension.
(Since the value for $v\notin\gMdy$ is given in \cite{dats} without a proof,
we choose to include the computation for any $v\in\gMf$.)

\section{The orbital volumes}\label{sec:vol}

Let $v$ be an arbitrary finite place.
In this section, we compute the volume
$\vol(\cK_vx)$ for standard representatives $x$ 
corresponding to quadratic ramified extensions.
Our purpose of this section is to prove the following.
\begin{prop}\label{vol}
Let $x\in V_{k_v}^\sst$ be one of the standard representatives
corresponding to quadratic ramified extensions. 
We have
\begin{equation*}
\vol(\cK_vx)=	2^{-1}|\Delta_{k_v(x)/k_v}|_v^{-1}(1-q_v^{-1})(1-q_v^{-2}).
\end{equation*}
\end{prop}

For the rest of this section we drop the subscript $v$
from $\cK_v,\co_v$ and write $\cK,\co$ instead.
We denote by $\gp$ the prime ideal of $\co$.
We consider the cases $v\notin\gMdy$ and
$v\in\gMdy$ simultaneously.
Let $\ord_v\colon k_v^\times\rightarrow \Z$
be the normalized discrete valuation.
We put $m=\ord_v(2)$. Note that $m=0$ if $v\notin\gMdy$.

Let $x=v_1^2+a_1v_1v_2+a_2v_2^2\in V_{k_v}^\sst$ be a standard representative
corresponding to a quadratic ramified extension.
In this case $x(v_1,1)$ is an Eisenstein polynomial
i.e. $a_1\in\gp$ and $a_2\in\gp\setminus\gp^2$.
We put $\Delta_{k_v(x)/k_v}=\gp^{\delta_x}$.
By definition,
$\delta_x=2\ord_v(a_1)$ if $1\leq \ord_v(a_1)\leq m$ and
$\delta_x=2m+1$ if $\ord_v(a_1)\geq m+1$.
Hence $\delta_x$ takes one of the values
$2, 4, \dots, 2m, 2m+1$.
(In the case $v\notin\gMdy$ and hence $m=0$, this should be counted as
$\delta_x$ only takes the value $1$.)
Let $x=\nr_{k_v(x)/k_v}(v_1+\varpi v_2)$.
Then $\varpi$ is a prime element of the ring of integers $\co_{k_x(x)}$
of $k_v(x)$.

Let $i$ be a positive integer.
For an $\co$-scheme $X$, let $\gr_{X,i}$ denote the reduction map
$X_{\co}\rightarrow X_{\co/\gp^i}$.
If the situation is obvious we drop $X$ and write $\gr_i$ instead.
For rational points $y_1,y_2\in X_{\co}$,
we use the notation $y_1\equiv y_2\;(\gp^i)$ if
$\gr_i(y_1)=\gr_i(y_2)$.
This coincides with the classical notation.

We put $n=\delta_x+2m+1$.
To compute $\vol(\cK x)$ we consider 
the congruence relation of modulo $\gp_v^n$
of the representation $(G_\co,V_\co)$.

\begin{defn}
We define $\cD_x=\{y\in V_{\co}\mid y\equiv x\,(\gp^n)\}$.
\end{defn}
\begin{lem}\label{Dcontainedrr}
We have $\cD_x\subset \cK x$.
\end{lem}
\begin{proof}
Let $y\in D_x$. First we show $y\in G_{k_v}x$.
Since $y\equiv x\;(\gp^n)$
and $\ord_v(P(x))=\delta_x$, we have $P(x)/P(y)\equiv 1\;(\gp^{2m+1})$
and so $P(y)/P(x)\in (k_v^\times)^2$.
Therefore the splitting fields of $x(v)$ and $y(v)$ coincide
and hence by Lemma \ref{rod}, we have $y\in G_{k_v}x$.
Let
$y=gx$, $g=(t,g_2)\in G_{k_v}$.
Note that
$|\chi(g)|_v=|t\det(g_2)|_v=1$
since $|P(x)|_v=|P(y)|_v$.
By Proposition \ref{Nstr},
multiplying an element of $N_{x\, k_v}$
if necessary, we may assume that
\begin{equation*}
|t|_v=|\det(g_2)|_v=1.
\end{equation*}
Then since $y(v)\in\sym^2\co^2$,
we have
\begin{align*}
y(1,0)&=t\n_{k_v(x)/k_v}(g_{211}+g_{212}\varpi)\in\co,\\
y(0,1)&=t\n_{k_v(x)/k_v}(g_{221}+g_{222}\varpi)\in\co.
\end{align*}
Hence both $\n_{k_v(x)/k_v}(g_{211}+g_{212}\varpi)$
and $\n_{k_v(x)/k_v}(g_{221}+g_{222}\varpi)$
are elements of $\co$ and so all entries of $g_2$ are in $\co$.
Since $|\det(g_2)|_v=1$, we conclude $g_2\in\gl(2)_{\co_v}$.
Hence $g\in \cK$ and the lemma follows.
\end{proof}

We will study the structure of $G_{x\,\co/\gp^n}$
following the method of \cite{kayuII}.
We first consider the structure of $N_x$.
\begin{prop}
We have $N_x\cong \co_{k_v(x)}^\times$ as a group scheme over $\co$.
\end{prop}
\begin{proof}
Let $R$ be any $\co$-algebra.
Then since $\co_{k_v(x)}=\co[\varpi]$,
we have $\co_{k_v(x)}\otimes_{\co} R=R[\varpi]$.
By a standard consideration,
we could see that the map
\begin{equation*}
\psi_{x\, R}\colon
N_{x\, R}\longrightarrow
R[\varpi]^\times,\qquad
n_x(c,d)
\longmapsto c+d\varpi
\end{equation*}
gives an isomorphism of these groups,
and this map satisfies the usual functorial property
with respect to homomorphism of $\co$-algebras.
This shows that there exists an isomorphism
$\psi_x\colon N_x\rightarrow \co_{k_v(x)}^\times$
as groups schemes over $\co$
such that $\psi_{x,R}$ is the induced isomorphism for all $R$.
\end{proof}

From this proposition we obtain the following.
\begin{cor}\label{ordNx}
The group $N_{x\,\co/\gp_v^n}$ is of order $q_v^{2n-1}(q_v-1)$.
\end{cor}

We now consider $G_{x\,\co/\gp^n}$.

\begin{lem}\label{stabrr}
We have $[G_{x\,\co/\gp^n}:N_{x\,\co/\gp^n}]=2q_v^{\delta_x}$.
\end{lem}
\begin{proof}
We shall count the number of elements of the right coset space
$N_{x\,\co/\gp^n}\backslash G_{x\,\co/\gp^n}$.
We claim that each right coset space contains
an element of the form
\begin{equation}\label{coset}
g=\left(1,\twtw 10us\right),
\qquad u\in\co/\gp^n,
\quad s\in(\co/\gp^n)^\times.
\end{equation}
Let $g\in G_{x\,\co/\gp^n}$.
Since $x(v)$ reduces to $v_1^2$ modulo $\gp$, 
we have $g_{221}\in\gp$
and therefore $g_{222}$ is a unit, which asserts
$A_x(g_{222},-g_{212})\in \gl(2)_{\co/\gp^n}$.
Hence multiplying 
$n_x(g_{222},-g_{212})\in N_x$ 
from the left if necessary,
we may assume $g_{212}=0$.
Then multiplying $n_x(g_{211}^{-1},0)$ if necessary
we further assume $g_{211}=1$.
Now comparing the coefficient of $v_1^2$ of $gx$ and $x$,
we have $t=1$. Therefore we have the claim.

Also it is easy to see that each coset has
exactly one representative in the form \eqref{coset}.
Hence we will consider when
such an element actually lies in $G_{x\,\co/\gp^n}$.
Suppose that $g$ is in the form \eqref{coset}.
By computation,
\begin{equation*}
(gx)(v_1,v_2)=v_1^2+(a_1s+2u)v_1v_2+(u^2+a_1us+a_2s^2)v_2^2.
\end{equation*}
Hence we can reduce the result of this lemma to the below. 
\end{proof}
\begin{lem}
Let $u\in\co/\gp^n$ and $s\in(\co/\gp^n)^\times$.
The system of the congruence equations
\begin{equation}\label{system}
a_1s+2u\equiv a_1\;(\gp^n),\quad
u^2+a_1us+a_2s^2\equiv a_2\;(\gp^n)
\end{equation}
with respect to $(u,s)$ has $2q_v^{\delta_x}$ numbers of solutions.
\end{lem}
\begin{proof}
In the analysis of this system it will be convenient to adopt
the usual abuse of notation by which classes in $\co/\gp^n$
and their representatives in $\co$ are denoted by the same symbol.
We put $\pi=a_2$, which is a uniformizer of $\co$ since $\ord_v(a_2)=1$.

We first consider the case $\delta_x=2m+1$,
hence $n=4m+2$.
In this case, we may assume $a_1=0$.
Then it is easy to see that the system \eqref{system}
is equivalent to
\begin{equation*}
u\equiv0\;(\gp^{3m+2}),\quad s^2\equiv1\;(\gp^{4m+1}).
\end{equation*}
Therefore for $u$ there are $q_v^m$ choices.
We consider the latter equation for $s$,
which is equivalent to
$\ord_v(s+1)+\ord_v(s-1)\geq 4m+1$.
Since $(s+1)-(s-1)=2$ and $\ord_v(2)=m$,
this condition is satisfied if and only if
$s\in\pm1+\gp^{3m+1}$.
Hence there are $2q_v^{m+1}$ possibilities for $s$.
Thus there are $2q_v^{2m+1}$ solutions of $(u,s)$ in all,
which proves the lemma in this case.

We next consider the case $\delta_x$ is one of $2,\dots,2m$,
which is occurred only if $v\in\gMdy$.
We put $l=\delta_x/2$, which is an integer between $1$ and $m$.
Note that $2/a_1\in\co$ since $\ord_v(a_1)=l\leq m$.
By definition, $n=2l+2m+1$.
From the second equation of \eqref{system}, 
we have $u\in\gp$.
We put $r=u/\pi$, which we regard as an element of $\co/\gp^{n-1}$.
Then the system \eqref{system} is equivalent to
\begin{equation*}
s+(2/a_1)\pi r\equiv 1\;(\gp^{l+2m+1}),\quad
\pi r^2+a_1sr+s^2=1\;(\gp^{2l+2m}).
\end{equation*}
We will consider the second equation under the first one.
From the first equation, we have
\begin{equation*}
a_1sr\equiv a_1r-2\pi r^2\;(\gp^{2l+2m+1})\quad\text{and}\quad
s^2\equiv 1-(4/a_1)\pi r+(4/a_1^2)\pi^2r^2\;(\gp^{l+3m+1})
\end{equation*}
and therefore by computation we see that
the second equation is equivalent to
\begin{equation*}
(4\pi^2/a_1^2-\pi)r^2+(a_1-4\pi/a_1)r\equiv0\;(\gp^{2l+2m}).
\end{equation*}
Let
$b_1=4\pi^2/a_1^2-\pi,
b_2=a_1-4\pi/a_1$, and 
$b=b_2/b_1$.
Then we have $\ord_v(b_1)=1,\ord_v(b_2)=l$ and $\ord_v(b)=l-1$,
and also the above equation is equivalent to
\begin{equation}\label{claim}
r(r+b)\equiv0\;(\gp^{2l+2m-1}).
\end{equation}
We claim that this equation is hold if and only if
$\ord_v(r)\geq l+2m$ or $\ord_v(r+b)\geq l+2m$.
In fact, if the equation is true,
either one order of $r$ and $r+b$ is greater than $l-1$.
Then since $\ord_v(b)=l-1$,
we have $\min\left\{\ord_v(r),\ord_v(r+b)\right\}=l-1$
and therefore $\max\left\{\ord_v(r),\ord_v(r+b)\right\}\geq l+2m$.
Conversely, under this condition we have
$\min\left\{\ord_v(r),\ord_v(r+b)\right\}=l-1$
and hence \eqref{claim} is hold.

All of these arguments shows that the system \eqref{system} is
hold if and only if either one condition of the following is valid:
\begin{itemize}
\item[(A)] $u\in \gp^{l+2m+1}$ and $s\in 1-2u/a_1+\gp^{l+2m+1}$,
\item[(B)] $u\in -b\pi+\gp^{l+2m+1}$ and $s\in 1-2u/a_1+\gp^{l+2m+1}$.
\end{itemize}
Since these conditions are mutually exclusive,
we have $2q_v^{2l}$ solutions in all and this finishes the proof.
\end{proof}

Now we give the proof of Proposition \ref{vol}
by putting together the result we have obtained before.

{\em Proof of Proposition \ref{vol}.}
Let $\gr_n$ be the reduction map $G_\co\rightarrow G_{\co/\gp^n}$.
Then by Lemma \ref{Dcontainedrr}, the set $\cK x=G_{\co}x$ is equal to
$\#(G_{\co}/\gr_n^{-1}(G_{x\,\co/\gp^n}))$ number of
disjoint copies of $\cD_x$.
Since
\begin{equation*}
G_{\co}/\gr_n^{-1}(G_{x\,\co/\gp^n})\cong
G_{\co/\gp^n}/G_{x\,\co/\gp^n},
\end{equation*}
by Lemma \ref{stabrr} and Corollary \ref{ordNx} we have
\begin{equation*}
\begin{aligned}
\vol(\cK x)
&=\vol(\cD_x)\cdot\frac{\#(G_{\co/\gp^n})}
{2q_v^{\delta_x}\cdot\#(N_{x\,\co/\gp^n})}\\
&=q_v^{-3n}\cdot
\frac{q_v^{n-1}(q_v-1)\cdot q_v^{4(n-1)}(q_v^2-1)(q_v^2-q_v)}
{2q_v^{\delta_x}\cdot q_v^{2n-1}(q_v-1)}\\
&=2^{-1}q_v^{-\delta_{x}}(1-q_v^{-1})(1-q_v^{-2}).
\end{aligned}
\end{equation*}
Since $|\Delta_{k(x)/k}|_v=q_v^{\delta_x}$, we obtained the desired result.
\hspace{\fill}$\square$
\begin{rem}
It is well known that there are
$2q_v^{l-1}(q_v-1)$
numbers of quadratic extensions of $k_v$
with the absolute value of the relative discriminant $q_v^{2l}$
for $1\leq l\leq m$ and
$2q_v^{m}$
numbers of quadratic extensions of $k_v$
with the absolute value of the relative discriminant $q_v^{2m+1}$.
Hence we have
\begin{equation*}
\begin{aligned}
\sum_{2\leq \del_{x,v}= 2l\leq 2m_v}\vol(\cK x) & = 
q_v^{-l}(1-q_v^{-1})^2(1-q_v^{-2}), \\
\sum_{\del_{x,v}= 2m+1}\vol(\cK x) & = 
q_v^{-(m+1)}(1-q_v^{-1})(1-q_v^{-2}),
\end{aligned}
\end{equation*}
where $x$ runs through all the standard representative
with the given condition of discriminants.
This result matches to \cite[Proposition 4.3]{dats}.
\end{rem}

\bibliographystyle{plain}

\begin{thebibliography}{1}

\bibitem{dats}
B.~Datskovsky.
\newblock A mean value theorem for class numbers of quadratic extensions.
\newblock {\em Contemporary Mathematics}, 143:179--242, 1993.

\bibitem{gauss}
C.F. Gauss.
\newblock {\em Disquisitiones arithmeticae}.
\newblock Yale University Press, New Haven, London, 1966.

\bibitem{goho}
D.~Goldfeld and J.~Hoffstein.
\newblock Eisenstein series of $1/2$-integral weight and the mean value of real
  {D}irichlet series.
\newblock {\em Invent. Math.}, 80:185--208, 1985.

\bibitem{kayuII}
A.C. Kable and A.~Yukie.
\newblock The mean value of the product of class numbers of paired quadratic
  fields, {II}.
\newblock {\em J. Math. Soc. Japan}, 55:739--764, 2003.

\bibitem{mufo}
D.~Mumford and J.~Fogarty.
\newblock {\em Geometric invariant theory}.
\newblock Springer-{V}erlag, Berlin, Heidelberg, New York, 2nd edition, 1982.

\bibitem{Shintanib}
T.~Shintani.
\newblock On zeta-functions associated with vector spaces of quadratic forms.
\newblock {\em J. Fac. Sci. Univ. Tokyo, Sect IA}, 22:25--66, 1975.

\bibitem{siegele}
C.L. Siegel.
\newblock The average measure of quadratic forms with given discriminant and
  signature.
\newblock {\em Ann. of Math.}, 45:667--685, 1944.

\bibitem{tanic}
T.~Taniguchi.
\newblock A mean value theorem for the square of class numbers of quadratic
  fields.
\newblock 2004.
\newblock in preparation.

\end{thebibliography}

\end{document}